\documentclass[reqno]{amsart}

\usepackage{amsmath}
\usepackage{amsfonts,amssymb}
\usepackage[dvips]{graphicx}
\usepackage{epsfig}
\newtheorem{theorem}{Theorem}[section]
\newtheorem{remark}{Remark}[section]

\newtheorem{definition}{Definition}[section]

\newtheorem{proposition}[theorem]{Proposition}

\numberwithin{equation}{section}
\begin{document}

\title[the compressible Euler equations in the isentropic nozzle flows ]
	{Remarks on the energy inequality of a global $L^{\infty}$ solution to the compressible Euler equations for the isentropic nozzle flow}
\author{Naoki Tsuge}
\address{Department of Mathematics Education, 
Faculty of Education, Gifu University, 1-1 Yanagido, Gifu
Gifu 501-1193 Japan.}
\email{tuge@gifu-u.ac.jp}
\thanks{
N. Tsuge's research is partially supported by Grant-in-Aid for Scientific 
Research (C) 17K05315, Japan.
}
\keywords{The Compressible Euler Equation, the nozzle flow, the compensated compactness, finite energy solutions, the
	modified Lax Friedrichs scheme.}
\subjclass{Primary 
35L03, 
35L65, 
35Q31, 
76N10,
76N15; 
Secondary
35A01, 
35B35,   
35B50, 
35L60,   
76H05,   
76M20.   
}
\date{}

\maketitle
\begin{abstract}
We study the compressible Euler equations in the isentropic nozzle flow. 
The global existence of an $L^{\infty}$ solution has been proved in (Tsuge in Nonlinear Anal. Real World Appl. 209: 217-238 (2017))
for large data and general nozzle. However, unfortunately, this solution does not
possess finiteness of energy. Although the modified Godunov scheme is introduced in 
this paper, we cannot deduce the energy inequality for the approximate solutions.

Therefore, our aim in the present paper is to derive the energy inequality for an $L^{\infty}$ solution. 
To do this, we introduce the modified Lax Friedrichs scheme, which has a recurrence formula consisting of 
discretized approximate solutions. We shall first deduce from the formula 
the energy inequality. Next, applying the compensated compactness method, the approximate solution converges to a weak solution. The energy inequality also holds for the solution as the limit. As a result, since our solutions are $L^{\infty}$, they possess finite energy and propagation, which are essential to physics.


\end{abstract}


\section{Introduction}
The present paper is concerned with isentropic gas flow in a nozzle. 
This motion is governed by the following compressible Euler equations:
\begin{equation}\left\{\begin{array}{ll}
\displaystyle{\rho_t+m_x=a(x)m,}\\
\displaystyle{m_t+\left(\frac{m^2}{\rho}+p(\rho)\right)_x
=a(x)\frac{m^2}{\rho},\quad x\in{\bf R}},
\end{array}\right.
\label{eqn:nozzle}
\end{equation}
where $\rho$, $m$ and $p$ are the density, the momentum and the 
pressure of the gas, respectively. If $\rho>0$, 
$v=m/\rho$ represents the velocity of the gas. For a barotropic gas, 
$p(\rho)=\rho^\gamma/\gamma$, where $\gamma\in(1,5/3]$ is the 
adiabatic exponent for usual gases. The given function $a(x)$ is 
represented by 
\begin{align*}
a(x)=-A'(x)/A(x)\quad\text{with}\quad A(x)=e^{-\int^x a(y)dy},
\end{align*}
where $A\in C^2({\bf R})$ is a slowly variable cross section area at $x$ in the nozzle.

We consider the Cauchy problem (\ref{eqn:nozzle}) with the initial data 
\begin{align}  
(\rho,m)|_{t=0}=(\rho_0(x),m_0(x)).
\label{eqn:I.C.}
\end{align}
The above problem \eqref{eqn:nozzle}--\eqref{eqn:I.C.} can be written in the following form 
\begin{align}\left\{\begin{array}{lll}
u_t+f(u)_x=g(x,u),\quad{x}\in{\bf R},\\
u|_{t=0}=u_0(x),
\label{eqn:IP}
\end{array}\right.
\end{align}
by using  $u={}^t(\rho,m)$, $\displaystyle f(u)={}^t\!\left(m, \frac{m^2}{\rho}+p(\rho)\right)$ and 
$\displaystyle{g(x,u)={}^t\!\left(a(x)m,a(x)\frac{m^2}{\rho}\right)}$.
The nozzle flow is applied in the various area, engineering, physics. Moreover, it is known that 
it is closely related to the flow of the solar wind. The detail can be found in \cite{T9}.

In the present paper, we consider an unsteady isentropic gas flow in particular. Let us survey the related mathematical results for 
the nozzle flow. The pioneer work in this direction is { Liu} \cite{L1}.  In 
\cite{L1}, {Liu} proved the existence of global solutions coupled with steady states, by the Glimm scheme, provided that the initial data have small total variation and are away from the sonic state.  
Recently, the existence theorems that include the transonic state have been obtained. 
The author \cite{T2} proved the global existence of solutions for 
the spherically symmetric case ($A(x)=x^2$ in \eqref{eqn:nozzle}) by the compensated compactness framework. {Lu} \cite{L4}, {Gu} and {Lu} \cite{LG} extended \cite{T2} to the nozzle flow with a monotone cross section area and the general pressure by using the vanishing viscosity method. In addition, the author \cite{T3} treated the Laval nozzle, which is a divergent and convergent nozzle. In these papers, the monotonicity of the cross section area plays an important role. For the general nozzle, the author \cite{T3} and \cite{T4} proved the global existence of a solution, provided that $a\in L^1({\bf R})$.

However, unfortunately, these solutions \cite{T2}, \cite{T3}, \cite{T4} and \cite{T8} do not possess finiteness of energy. Although the modified Godunov scheme is introduced in these paper, we cannot deduce the energy inequality for the corresponding approximate solution. 
Since our solutions are weak ones, which are defined almost everywhere, it is difficult to deduce the energy inequality for the weak solutions directly. Our main purpose of the present paper is to prove the inequality for solutions. Our strategy is as follows. 
We introduce the modified Lax Friedrichs scheme. By using the scheme, we can obtain the global existence of a solution 
in a similar manner to the modified Godunov scheme. Moreover, this has a recurrence formula consisting of discretized approximate solutions (see \eqref{eqn:energy1}). 
We shall first deduce from the formula the energy inequality. Since it consists of discretized values such as sequence, the treatment is comparetively easy.  
Next, applying the compensated compactness, the approximate 
solutions converge to a weak solution. As a result, the energy inequality also holds for the weak solution as the limit. This idea is employed in \cite{T10} and \cite{T11}. 
In this paper, we prove the 
energy inequality for \cite{T8} in particular. However, we can similarly apply our method to the other cases \cite{T2}, \cite{T3} and \cite{T4}.

The above finite energy solutions have recently received attention in \cite{CS} and \cite{LW}. In these results, solutions are constructed in $L^{p}$. On the other hand, our solution is $L^{\infty}$, which yields finite propagation. 
Therefore, our solution possesses finiteness of both energy and propagation, which are essential to physics.

To state our main theorem, we define the Riemann invariants $w,z$, which play important roles
in this paper, as
\begin{definition}
\begin{align*}
w:=\frac{m}{\rho}+\frac{\rho^{\theta}}{\theta}=v+\frac{\rho^{\theta}}{\theta},
\quad{z}:=\frac{m}{\rho}-\frac{\rho^{\theta}}{\theta}
=v-\frac{\rho^{\theta}}{\theta}\quad(\theta:=(\gamma-1)/2).
\end{align*}
\end{definition}
These Riemann invariants satisfy the following.
\begin{remark}\label{rem:Riemann-invariant}
\normalfont
\begin{align}
&|w|\geq|z|,\;w\geq0,\;\mbox{\rm when}\;v\geq0.\quad
|w|\leq|z|,\;z\leq0,\;\mbox{\rm when}\;v\leq0.\label{eqn:inequality-Riemann}\\
&v=\frac{w+z}2,
\;\rho=\left(\frac{\theta(w-z)}2\right)^{1/\theta},\;m=\rho v.
\label{eqn:relation-Riemann}
\end{align}From the above, the lower bound of $z$ and the upper bound of $w$ yield the bound of $\rho$ and $|v|$.
\end{remark}

Moreover, we define the entropy weak solution.
\begin{definition}
A measurable function $u(x,t)$ is called a global {\it entropy weak solution} of the 
Cauchy problems \eqref{eqn:IP} if 
\begin{align*}
\int^{\infty}_{-\infty}\int^{\infty}_0u\phi_t+f(u)\phi_x+g(x,u)\phi dxdt
+\int^{\infty}_{-\infty}u_0(x)\phi(x,0)dx=0
\end{align*}
holds for any test function $\phi\in C^1_0({\bf R}\times{\bf R}_+)$ and 
\begin{align}
\int^{\infty}_{-\infty}\int^{\infty}_0\hspace{-1ex}\eta(u)\psi_t+q(u)\psi_x+\nabla\eta(u) g(x,u)\psi dxdt+\int^{\infty}_{-\infty}\hspace{-1ex}\eta(u_0(x))\psi(x,0)dx
\geq0
\label{eqn:entropy inequality}
\end{align}
holds for any non-negative test function $\psi\in C^1_0({\bf R}\times{\bf R}_+)$, where 
$(\eta,q)$ is a pair of convex entropy--entropy flux of \eqref{eqn:nozzle} (see Section 2).
\end{definition}

We assume the following.\\
There exists  a nonnegative function $b\in C^1({\bf R})$ such that  
\begin{align}
\begin{split}
|a(x)|\leq \mu b(x),\quad
\max\left\{\int^{\infty}_0b(x)dx,\;\int^0_{-\infty}b(x)dx\right\}\leq\frac12\log\frac1{\sigma},
\end{split}
\label{eqn:condition-M}
\end{align}
where $\mu=\frac{(1-\theta)^2}{\theta(1+\theta-2\sqrt{\theta})},\;\sigma
=\frac{1-\theta}{(1-\sqrt{\theta})(2\sqrt{\theta+1}+\sqrt{\theta}-1)}$. Here we notice that 
$0<\sigma<1$.


From the similar argument of \cite{T8}, we have
\begin{theorem}\label{thm:main}
We assume that, for $b$ in \eqref{eqn:condition-M}  and any fixed nonnegative constant $M$, 
initial density and momentum data 
$u_0=({\rho}_0, {m}_0)\in{L}^{\infty}({\bf R})$ satisfy
\begin{align}
0\leq\rho_0(x)
,\;\; -Me^{-\int^x_0b(y)dy}\leq{z}(u_0(x)),\;\; w(u_0(x))\leq Me^{\int^x_0b(y)dy}\quad\text{\rm a.\,e. }x\in{\bf R}.
\label{eqn:IC}
\end{align}

Then the Cauchy problem $(\ref{eqn:IP})$ has a 
global entropy weak solution $u(x,t)$ satisfying  the same inequalities as \eqref{eqn:IC}
\begin{align*}
&0\leq\rho(x,t)
,\;\; -Me^{-\int^x_0b(y)dy}\leq{z}(u(x,t)),\;\; w(u(x,t))\leq Me^{\int^x_0b(y)dy}\\
&\text{\rm a.\,e. }(x,t)
\in{\bf R}\times{\bf R}_+.
\end{align*}
\end{theorem}
\begin{remark}
In view of $\eqref{eqn:condition-M}_2$, \eqref{eqn:IC} implies that we can supply arbitrary $L^{\infty}$ data.
\end{remark}

Then, our main theorem is as follows.
\begin{theorem}\label{thm:main2}
If the energy of initial data $\int_{\bf R}A(x)\eta_*(u_0(x))dx$ is finite, for the solution of Theorem \ref{thm:main}, the following energy inequality holds.
\begin{align}
\int_{\bf R}A(x)\eta_*(u(x,t))dx\leq\int_{\bf R}A(x)\eta_*(u_0(x))dx\quad\text{\rm a.\,e. }t>0,
\label{eqn:energy}
\end{align}where $\eta^*$ is the mechanical energy defined as follows.
\begin{align*}
\eta_*:=\frac12\frac{m^2}{\rho}+\frac1{\gamma(\gamma-1)}\rho^{\gamma}
\end{align*}
\end{theorem}

The present paper is organized as follows.
In Section 2, we review the Riemann problem and the properties of Riemann 
solutions.  In Section 3, we construct approximate solutions by 
the modified Lax Friedrichs scheme. In Section 4, we drive the recurrence formula consisting of discretized approximate solutions. 
We shall deduce the energy inequality for the formula.

\section{Preliminary}
In this section, we first review some results of the Riemann solutions 
for the homogeneous system of gas dynamics. Consider the homogeneous system 
\begin{equation}\left\{\begin{array}{ll}
\rho_t+m_x=0,\\
\displaystyle{m_t+\left(\frac{m^2}{\rho}+p(\rho)\right)_x=0,
\quad{p}(\rho)=\rho^{\gamma}/\gamma.}
\end{array}\right.
\label{eqn:homogeneous}
\end{equation}

A pair of functions $(\eta,q):{\bf R}^2\rightarrow{\bf R}^2$ is called an 
entropy--entropy flux pair if it satisfies an identity
\begin{equation}
\nabla{q}=\nabla\eta\nabla{f}.
\label{eqn:eta-q}
\end{equation}
Furthermore, if, for any fixed ${m}/{\rho}\in(-\infty,\infty)$, $\eta$ 
vanishes on the vacuum $\rho=0$, then $\eta$ is called a {\it weak entropy}. 
For example, the mechanical energy--energy flux pair 
\begin{equation}
\eta_*:=\frac12\frac{m^2}{\rho}+\frac1{\gamma(\gamma-1)}\rho^{\gamma},\quad
q_*:=m\left(\frac12\frac{m^2}{\rho^2}+\frac{\rho^{\gamma-1}}{\gamma-1}\right)  
\label{eqn:mechanical}
\end{equation}
should be a strictly convex weak entropy--entropy flux pair. 

The jump discontinuity in a weak solutions to (\ref{eqn:homogeneous}) must satisfy the following Rankine--Hugoniot condition 
\begin{align}
\lambda(u-u_0)=f(u)-f(u_0),
\label{eqn:R-H}
\end{align}
where $\lambda$ is the propagation speed of the discontinuity, 
$u_0=(\rho_0,m_0)$ and $u=(\rho,m)$ are the corresponding left and 
right state, respectively. 
A jump discontinuity is called a {\it shock} if it satisfies the entropy 
condition 
\begin{align}
\lambda(\eta(u)-\eta(u_0))-(q(u)-q(u_0))\geq0
\label{eqn:entropy-condition}
\end{align}
for any convex entropy pair $(\eta,q)$.

There are two distinct types of rarefaction and shock curves in the 
isentropic gases.
Given a left state $(\rho_0,m_0)$ or $(\rho_0,v_0)$, the possible states 
$(\rho,m)$ or $(\rho,v)$ that can be connected to $(\rho_0,m_0)$ or 
$(\rho_0,v_0)$ on the right by a rarefaction or a shock curve form
a 1-rarefaction wave curve $R_1(u_0)$, a 2-rarefaction wave curve $R_2(u_0)$,
a 1-shock curve $S_1(u_0)$ and a 2-shock curve $S_2(u_0)$: 
\begin{align*}
&R_1(u_0):w=w_0,\;\rho<\rho_0,\quad R_2(u_0):z=z_0,\;\rho>\rho_0,\\
&S_1(u_0):
\displaystyle{v-v_0=-\sqrt{\frac1{\rho\rho_0}\frac{p(\rho)-p(\rho_0)}
{\rho-\rho_0}}(\rho-\rho_0)\quad\rho>\rho_0>0,}\\
&S_2(u_0):
\displaystyle{v-v_0=\sqrt{\frac1{\rho\rho_0}\frac{p(\rho)-p(\rho_0)}
{\rho-\rho_0}}(\rho-\rho_0)\quad\rho<\rho_0,}
\end{align*}
respectively. Here we notice that shock wave curves are deduced from the\linebreak 
Rankine--Hugoniot condition (\ref{eqn:R-H}).

\begin{figure}[htbp]
\begin{center}
\vspace{-2ex}
\hspace{-6ex}
\includegraphics[scale=0.42]{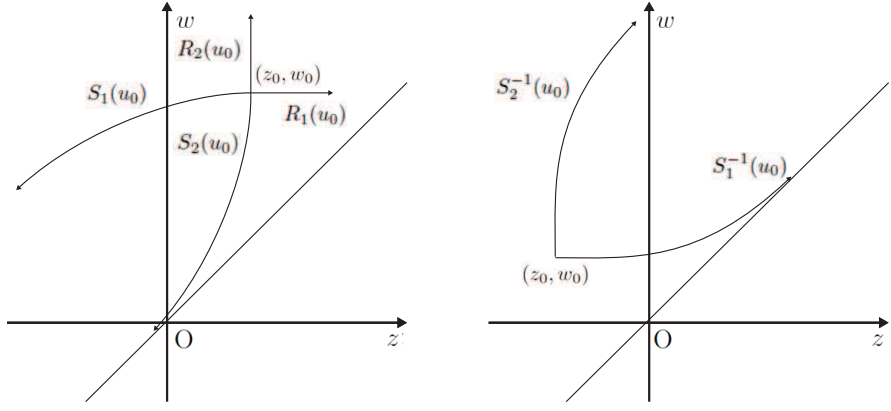}
\end{center}
\caption{The rarefaction curves, the shock curves and the inverse rarefaction curves in $(z,w)$-plane}
\end{figure}


\subsection{Riemann Solution}
Given a right state $(\rho_0,m_0)$ or $(\rho_0,v_0)$, the possible states 
$(\rho,m)$ or $(\rho,v)$ that can be connected to $(\rho_0,m_0)$ or 
$(\rho_0,v_0)$ on the left by a shock curve constitute
1-inverse shock curve $S_1^{-1}(u_0)$ and 2-inverse shock curve 
\begin{align*}&&S_1^{-1}(u_0):
\displaystyle{v-v_0=-\sqrt{\frac1{\rho\rho_0}\frac{p(\rho)-p(\rho_0)}
{\rho-\rho_0}}(\rho-\rho_0),\quad\underline{\rho<\rho_0},}\\
&&S_2^{-1}(u_0):
\displaystyle{v-v_0=\sqrt{\frac1{\rho\rho_0}\frac{p(\rho)-p(\rho_0)}
{\rho-\rho_0}}(\rho-\rho_0),\quad\underline{\rho>\rho_0>0},}
\end{align*}
respectively.

Next we define a rarefaction shock. Given 
$u_0,u$ on $S_i^{-1}(u_0)\;(i=1,2)$, 
we call the piecewise 
constant solution to (\ref{eqn:homogeneous}), which 
consists of the left and right states $u_0,u$ a {\it rarefaction shock}. 
Here, notice the following: although the inverse shock curve has the same 
form as the shock curve, the underline expression in $S_i^{-1}(u_0)$ 
is different from the corresponding part in $S_i(u_0)$.
Therefore the rarefaction shock does not satisfy the entropy condition.

We shall use a rarefaction shock in approximating a rarefaction wave.
In particular, when we consider a rarefaction shock, we call the inverse shock 
curve connecting $u_0$ and $u$ a {\it rarefaction shock curve}.

\begin{figure}[htbp]
	\begin{center}
		\vspace{-2ex}
		\hspace{-6ex}
		\includegraphics[scale=0.42]{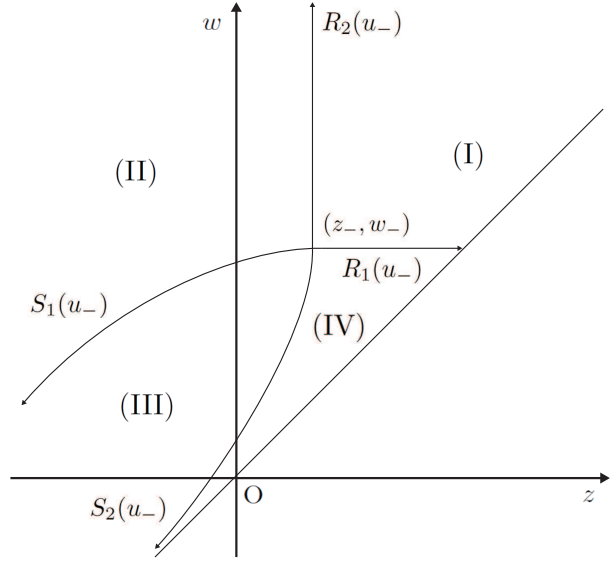}
	\end{center}
	\caption{The elementary wave curves in $(z,w)$-plane}
	\label{Fig:Riemann}
\end{figure}

From the properties of these curves in phase plane $(z,w)$, we can construct 
a unique solution for the Riemann problem 
\begin{equation}
u|_{t=0}=\left\{\begin{array}{ll}
u_-,\quad{x}<x_0,\\
u_+,\quad{x}>x_0,
\label{eqn:Riemann}
\end{array}\right.\end{equation} 
where $x_0\in(-\infty,\infty)$, $\rho_{\pm}\geq0$ and $m_{\pm}$ are constants 
satisfying $|m_{\pm}|\leq{C}\rho_{\pm}$. The Riemann solution consists of the following (see Fig. \ref{Fig:Riemann}).
\begin{enumerate}
\item $(z_+,w_+)\in$ (I): 1-rarefaction curve and 2-rarefaction curve; 
\item $(z_+,w_+)\in$ (II): 1-shock curve and 2-rarefaction curve; 
\item $(z_+,w_+)\in$ (III): 1-shock curve and 2-shock curve; 
\item $(z_+,w_+)\in$ (IV): 1-rarefaction curve and 2-shock curve,
\end{enumerate}
where $z_{\pm}=m_{\pm}/\rho_{\pm}-(\rho_{\pm})^{\theta}/\theta,\;w_{\pm}=m_{\pm}/\rho_{\pm}+(\rho_{\pm})^{\theta}/\theta$ respectively.

We denote the solution the Riemann solution $(u_-,u_+)$.

\section{Construction of Approximate Solutions}
\label{sec:construction-approximate-solutions}
In this section, we construct approximate solutions. In the strip 
$0\leq{t}\leq{T}$ for any fixed $T\in(0,\infty)$, we denote these 
approximate solutions by $u^{\varDelta}(x,t)
=(\rho^{\varDelta}(x,t),m^{\varDelta}(x,t))$. 
Let ${\varDelta}x$ and ${\varDelta}{t}$ be the space 
and time mesh lengths, respectively. Moreover, for any fixed positive 
value $X$, we assume that 
\begin{align}
A(x) \text{ is a constant in } |x|>X.
\label{eqn:X}
\end{align} 
Then we notice that $a(x)$ is bounded and has a compact support.

Let us define the approximate solutions by using the modified Lax Friedrichs scheme. 
We set 
\begin{align*}
(j,n)\in\mbox{\bf Z}\times\mbox{\bf Z}_{\geq0}.
\end{align*}
In addition, using $M$ in \eqref{eqn:IC}, we take ${\varDelta}x$ and ${\varDelta}{t}$ such that 
\begin{align*}
\frac{{\varDelta}x}{{\varDelta}{t}}=2Me^{\max\left\{\int^{\infty}_0b(x)dx,\;\int^0_{-\infty}b(x)dx\right\}}.
\end{align*}

First we define $u^{\varDelta}(x,-0)$ by 
\begin{align*}
u^{\varDelta}(x,-0)=u_0(x)
\end{align*}\begin{eqnarray*}
	u^{\varDelta}(x,-0)=\chi_{\scriptscriptstyle{X}}(x)u_0(x),
\end{eqnarray*}
where 
\begin{equation}
\begin{split}
\chi_{\scriptscriptstyle{X}}(x)=\left\{
\begin{array}{ll}
1,\quad x\leq X,\\
0,\quad x\geq X
\end{array}
\right.
\end{split}
\label{eqn:cutoff}
\end{equation}
and set 
\begin{align*}
J_n=\{k;k+n=\text{even}\}.
\end{align*}
Then, for $j\in J_0$, we define $E_j^0(u)$ by
\begin{align*}
E_j^0(u)=\frac1{2{\varDelta}x}\int_{{(j-1)}{\varDelta}x}^{(j+1){\varDelta}x}
u^{\varDelta}(x,-0)dx.
\end{align*}

Next, assume that $u^{\varDelta}(x,t)$ is defined for $t<n{\varDelta}{t}$. 
Then, for $j\in J_n$, we define $E^n_j(u)$ by 
\begin{align*}
E^n_j(u)=\frac1{2{\varDelta}x}\int_{{(j-1)}{\varDelta}x}^{(j+1){\varDelta}x}u^{\varDelta}(x,n{\varDelta}{t}-0)dx.
\end{align*}

Moreover, for $j\in J_n$, we define $u_j^n=(\rho_j^n,m_j^n)$ as follows.\\
We choose $\delta$ such that $1<\delta<1/(2\theta)$. If 
\begin{align*}
E^n_j(\rho):=
\frac1{2{\varDelta}x}\int_{{(j-1)}{\varDelta}x}^{(j+1){\varDelta}x}\rho^{\varDelta}(x,n{\varDelta}{t}-0)dx<({\varDelta}x)^{\delta},
\end{align*} 
we define $u_j^n$ by $u_j^n=(0,0)$;
otherwise, setting

\begin{align}
\begin{split}
{z}_j^n:&=
\max\left\{z(E_j^n(u)),\;-Me^{-\int^{j{\varDelta}x}_0b(x)dx}\right\}
\\&\hspace*{100pt}\mbox{and}\\
w_j^n:&=\min\left\{w(E_j^n(u)),\;Me^{\int^{j{\varDelta}x}_0b(x)dx}\right\}
,
\end{split}
\label{eqn:def-u^n_j}
\end{align} 
we define $u_j^n$ by
\begin{align*}
u_j^n:=(\rho_j^n,m_j^n)
:=\left(\left\{\frac{\theta(w_j^n-z_j^n)}{2}\right\}
^{1/\theta},
\left\{\frac{\theta(w_j^n-z_j^n)}{2}\right\}^{1/\theta}
\frac{w_j^n+z_j^n}{2}\right).
\end{align*}

\subsection{Construction of Approximate Solutions in the Cell}
\label{subsec:construction-approximate-solutions}
By using $u_j^n$ defined above, we 
construct the approximate solutions with 
$u^{\varDelta}((j-1){\varDelta}x,n{\varDelta}t+0)=u_{j-1}^n$ and $u^{\varDelta}((j+1){\varDelta}x,n{\varDelta}t+0)=u_{j+1}^n$in the cell $(j-1){\varDelta}x\leq{x}<(j+1){\varDelta}x,\;n{\varDelta}{t}\leq{t}<(n+1){\varDelta}{t}\quad(j\in J_{n+1},\;n\in{\bf Z}_{\geq0})$.

We first solve a Riemann problem with initial data $(u_{j-1}^n,u_{j+1}^n)$. 
Call constants $u_{\rm L}(=u_{j-1}^n), u_{\rm M}, u_{\rm R}(=u_{j+1}^n)$ the left, middle and 
right states, respectively. Then the following four cases occur.
\begin{itemize}
\item {\bf Case 1} A 1-rarefaction wave and a 2-shock arise. 
\item {\bf Case 2} A 1-shock and a 2-rarefaction wave arise. 
\item {\bf Case 3} A 1-rarefaction wave and a 2-rarefaction arise.
\item {\bf Case 4} A 1-shock and a 2-shock arise.
\end{itemize}
We then construct approximate solutions $u^{\varDelta}(x,t)$ by perturbing 
the above Riemann solutions. We consider only the case in which $u_{\rm M}$ is away from the vacuum. The other case (i.e., the case where $u_{\rm M}$ is near the vacuum) is a little technical. Therefore, we postpone the case near the vacuum to Appendix A. 

\vspace*{10pt}
{\bf The case where $u_{\rm M}$ is away from the vacuum}

Let $\alpha$ be a constant satisfying $1/2<\alpha<1$. Then we can choose 
a positive value $\beta$ small enough such that $\beta<\alpha$, $1/2+\beta/2<\alpha<
1-2\beta$, $\beta<2/(\gamma+5)$ and $(9-3\gamma)\beta/2<\alpha$.

We first consider the case where $\rho_{\rm M}>({\varDelta}x)^{\beta}$, 
which  means $u_{\rm M}$ is away from the vacuum. In this step, we 
consider Case 1 in particular. The constructions of Cases 2--4 are similar 
to that of Case 1.

Consider the case where a 1-rarefaction wave and a 2-shock arise as a Riemann 
solution with initial data $(u_{j-1}^n,u_{j+1}^n)$. Assume that 
$u_{\rm L},u_{\rm M}$ 
and $u_{\rm M},u_{\rm R}$ are connected by a 1-rarefaction and a 2-shock 
curve, respectively. 

{\it Step 1}.\\
In order to approximate a 1-rarefaction wave by a piecewise 
constant {\it rarefaction fan}, we introduce the integer  
\begin{align*}
p:=\max\left\{[\hspace{-1.2pt}[(z_{\rm M}-z_{\rm L})/({\varDelta}x)^{\alpha}]
\hspace{-1pt}]+1,2\right\},
\end{align*}
where $z_{\rm L}=z(u_{\rm L}),z_{\rm M}=z(u_{\rm M})$ and $[\hspace{-1.2pt}[x]\hspace{-1pt}]$ is the greatest integer 
not greater than $x$. Notice that
\begin{align}
p=O(({\varDelta}x)^{-\alpha}).
\label{eqn:order-p}
\end{align}
Define \begin{align*}
z_1^*:=z_{\rm L},\;z_p^*:=z_{\rm M},\;w_i^*:=w_{\rm L}\;(i=1,\ldots,p),
\end{align*}
and 
\begin{align*}
z_i^*:=z_{\rm L}+(i-1)({\varDelta}x)^{\alpha}\;(i=1,\ldots,p-1).
\end{align*}
We next introduce the rays $x=j{\varDelta}x+\lambda_1(z_i^*,z_{i+1}^*,w_{\rm L})
(t-n{\varDelta}{t})$ separating finite constant states 
$(z_i^*,w_i^*)\;(i=1,\ldots,p)$, 
where  
\begin{align*}
\lambda_1(z_i^*,z_{i+1}^*,w_{\rm L}):=v(z_i^*,w_{\rm L})
-S(\rho(z_{i+1}^*,w_{\rm L}),\rho(z_i^*,w_{\rm L})),
\end{align*}
\begin{align*}
\rho_i^*:=\rho(z_i^*,w_{\rm L}):=\left(\frac{\theta(w_{\rm L}-z_i^*)}2\right)^{1/\theta}\;,
\quad{v}_i^*:={v}(z_i^*,w_{\rm L}):=\frac{w_{\rm L}+z_i^*}2
\end{align*}
and

\begin{align}
S(\rho,\rho_0):=\left\{\begin{array}{lll}
\sqrt{\displaystyle{\frac{\rho(p(\rho)-p(\rho_0))}{\rho_0(\rho-\rho_0)}}}
,\quad\mbox{if}\;\rho\ne\rho_0,\\
\sqrt{p'(\rho_0)},\quad\mbox{if}\;\rho=\rho_0.
\end{array}\right.
\label{eqn:s(,)}
\end{align}

We call this approximated 1-rarefaction wave a {\it 1-rarefaction fan}.

\vspace*{10pt}
{\it Step 2}.\\
In this step, we replace the above constant states  with the following functions of $x$:

\begin{definition}\label{def:steady-state}\normalfont
	For given constants $x_d$ satisfying  $(j-1){\varDelta}x\leq{x_d}\leq(j+1){\varDelta}x$ and 
	\begin{align}
	\begin{split}
	(z_d,w_d):=\left(\frac{m_d}{\rho_d}-\frac{(\rho_d)^{\theta}}
	{\theta},\frac{m_d}{\rho_d}+\frac{(\rho_d)^{\theta}}{\theta}\right)
	\quad\text{or}\quad
	u_d=(\rho_d,m_d)
	\end{split}
	\label{eqn:steady-state-data}
	\end{align}
	satisfying $|m_d|\leq{C}\rho_d$, we set 
	\begin{align*}
	z(x)=z_de^{-\int^x_{x_d}b(y)dy},\quad
	w(x)=w_de^{\int^x_{x_d}b(y)dy}.
	\end{align*}
	Using $z(x)$ and $w(x)$, we define  
	\begin{align}
	u(x)=(\rho(x),m(x))
	\label{eqn:solution-steady-state}
	\end{align}
	by the relation 
	\eqref{eqn:relation-Riemann} as follows.
	\begin{align}
	v(x)=\frac{w(x)+z(x)}2,\;\rho(x)=\left(\frac{\theta(w(x)-z(x))}2\right)^{1/\theta},\;m(x)=\rho(x) v(x).
\label{eqn:relation-Riemann2}
	\end{align}
We then define $\bar{\mathcal U}(x,x_d,u_d)$ with data $u_d$ at $x_d$ as  
\eqref{eqn:solution-steady-state} .
\end{definition}

Moreover, for given functions $\bar{u}(x)$, we define $z(x,t)$ and $w(x,t)$ by 
\begin{align*}
\begin{split}
z(x,t)=&\bar{z}(x)
-\left\{a(x)\bar{v}(x)(\bar{\rho}(x))^{\theta}
-b(x)\lambda_1(\bar{u}(x))\bar{z}(x)\right\}
(t-n{\varDelta{t}}),\\
w(x,t)=&\bar{w}(x)+
\left\{a(x)\bar{v}(x)(\bar{\rho}(x))^{\theta}
-b(x)\lambda_2(\bar{u}(x))\bar{w}(x)\right\}(t-n{\varDelta{t}}).
\end{split}
\end{align*}
Then, using $z(x,t)$ and $w(x,t)$, we define $u(x,t)=(\rho(x,t),m(x,t))$ in a similar manner to
\eqref{eqn:relation-Riemann2}.  We denote $u(x,t)$ by ${\mathcal U}(x,t;\bar{u})$

Let $\bar{u}_{\rm L}(x)$ and $\bar{u}_{\rm R}(x)$ be $\bar{\mathcal U}
(x,(j-1){\varDelta}x,u_{\rm L})$ and 
$\bar{\mathcal U}(x,(j+1){\varDelta}x,u_{\rm R})$, respectively. Set 
$
\bar{u}_1(x):=\bar{u}_{\rm L}(x),\;u_1(x,t)={\mathcal U}(x,t;\bar{u}_1),\;u_{\rm R}(x,t)={\mathcal U}(x,t;\bar{u}_{\rm R})$ and $x_1:=(j-1){\varDelta}x.
$

\begin{figure}[htbp]
\begin{center}
\hspace{-2ex}
\includegraphics[scale=0.3]{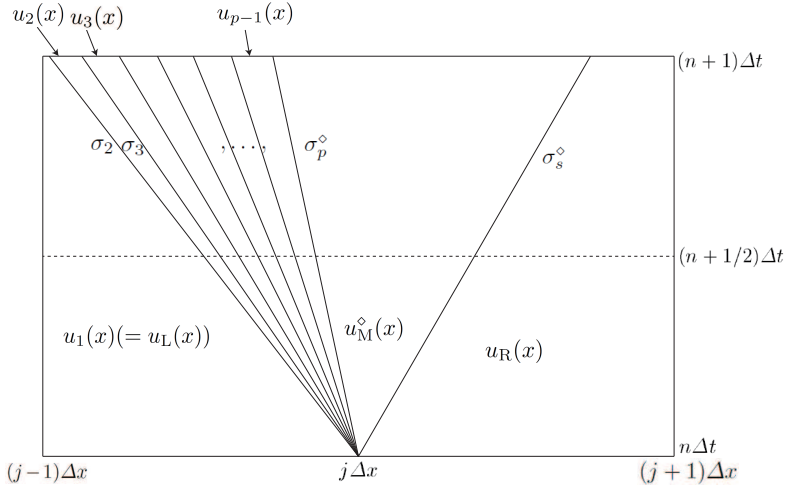}
\end{center}
\caption{The approximate solution in the case where a 1-rarefaction and 
a 2-shock arise in the cell.}
\label{case1-1cell}
\end{figure}

First, by the implicit function theorem, we determine a propagation speed $\sigma_2$ and $u_2=(\rho_2,m_2)$ such that 
\begin{itemize}
	\item[(1.a)] $z_2:=z(u_2)=z^*_2$
	\item[(1.b)] the speed $\sigma_2$, the left state ${u}_1(x_2,(n+1/2){\varDelta}t)$ and the right state $u_2$ satisfy the Rankine--Hugoniot conditions, i.e.,
	\begin{align*}
	f(u_2)-f({u}_1(x_2,(n+1/2){\varDelta}t))=\sigma_2(u_2-{u}_1(x_2,(n+1/2){\varDelta}t)),
	\end{align*}
\end{itemize}
where $x_2:=j{\varDelta}x+\sigma_2{\varDelta}t/2$. Then we fill up by ${u}_1(x)$ the sector where $n{\varDelta}t\leq{t}<
(n+1){\varDelta}t,(j-1){\varDelta}x\leq{x}<{j}{\varDelta}x+
\sigma_2(t-n{\varDelta}t)$ (see Figure \ref{case1-1cell}) 
and set $\bar{u}_2(x)=\bar{\mathcal U}(x,x_2,u_2)$ and $u_2(x,t)={\mathcal U}(x,t;\bar{u}_2)$.

Assume that $u_k$, ${u}_k(x,t)$ and a propagation speed $\sigma_k$ with
$\sigma_{k-1}<\sigma_k$ are defined. Then we similarly determine
$\sigma_{k+1}$ and $u_{k+1}=(\rho_{k+1},m_{k+1})$ such that 
\begin{itemize}
	\item[($k$.a)] $z_{k+1}:=z(u_{k+1})=z^*_{k+1}$,
	\item[($k$.b)] $\sigma_{k}<\sigma_{k+1}$,
	\item[($k$.c)] the speed 
	$\sigma_{k+1}$, 
	the left state ${u}_k(x_{k+1},(n+1/2){\varDelta}t)$ and the right state $u_{k+1}$ satisfy 
	the Rankine--Hugoniot conditions, 
\end{itemize}
where 
$x_{k+1}:=j{\varDelta}x+\sigma_{k+1}{\varDelta}t/2$. Then we fill up by ${u}_k(x,t)$ the sector where
$n{\varDelta}t\leq{t}<(n+1){\varDelta}t,{j}{\varDelta}x+\sigma_k(t-{\varDelta}t)\leq{x}<{j}{\varDelta}x+\sigma_{k+1}(t-n{\varDelta}t)$  (see Figure \ref{case1-1cell}) and 
set $\bar{u}_{k+1}(x)=\bar{\mathcal U}(x,x_{k+1},u_{k+1})$ and $u_{k+1}(x,t)={\mathcal U}(x,t;\bar{u}_{k+1})$.  
By induction, we define $u_i$, ${u}_i(x,t)$ and $\sigma_i$ $(i=1,\ldots,p-1)$.
Finally, we determine a propagation speed $\sigma_p$ and $u_p=(\rho_p,m_p)$ such that
\begin{itemize}
	\item[($p$.a)] $z_p:=z(u_p)=z^*_p$,
	\item[($p$.b)] the speed $\sigma_p$, 
	and the left state ${u}_{p-1}(x_p,(n+1/2){\varDelta}t)$ and the right state $u_p$ satisfy the Rankine--Hugoniot conditions, 
	\end{itemize}where $x_p:=j{\varDelta}x+\sigma_p{\varDelta}t/2$. 
We then fill up by ${u}_{p-1}(x,t)$ and $u_p$ the sector where
$n{\varDelta}t\leq{t}<(n+1){\varDelta}t,{j}{\varDelta}x+\sigma_{p-1}
(t-n{\varDelta}t)\leq{x}<{j}{\varDelta}x+\sigma_{p}(t-n{\varDelta}t)$ 
and the line $n{\varDelta}t\leq{t}<(n+1){\varDelta}t,x={j}{\varDelta}x+\sigma_{p}(t-n{\varDelta}t)$, respectively.

Given $u_{\rm L}$ and $z_{\rm M}$ with $z_{\rm L}\leq{z}_{\rm M}$, we denote 
this piecewise functions of $x$ 1-rarefaction wave by 
$R_1^{\varDelta}(z_{\rm M})(u_{\rm L})$. Notice that from the construction 
$R^{\varDelta}_1(z_{\rm M})(u_{\rm L})$ connects $u_{\rm L}$ and $u_p$ 
with $z_p=z_{\rm M}$. 

Now we fix ${u}_{\rm R}(x,t)$ and ${u}_{p-1}(x,t)$. 
Let $\sigma_s$ be 
the propagation speed of the 2-shock connecting $u_{\rm M}$ and $u_{\rm R}$.
Choosing ${\sigma}^{\diamond}_p$ near to $\sigma_p$, ${\sigma}^{\diamond}_s$ 
near to 
$\sigma_s$ and $u^{\diamond}_{\rm M}$ near to $u_{\rm M}$, we fill up by ${u}^{\diamond}_{\rm M}(x,t)=
{\mathcal U}(x,t,\bar{u}^{\diamond}_{\rm M})$ the gap between $x=j{\varDelta}x+{\sigma}^{\diamond}_{p}
(t-n{\varDelta}{t})$ and $x=j{\varDelta}x+{\sigma}^{\diamond}_s(t-n{\varDelta}{t})$, such that 
\begin{itemize}
	\item[(M.a)] $\sigma_{p-1}<\sigma^{\diamond}_p<\sigma^{\diamond}_s$, 
	\item[(M.b)] the speed ${\sigma}^{\diamond}_p$, the left and right states 
	${u}_{p-1}(x^{\diamond}_{p},(n+1/2){\varDelta}t),{u}^{\diamond}_{\rm M}(x^{\diamond}_{p},(n+1/2){\varDelta}t)$ 
	satisfy the Rankine--Hugoniot conditions,
	\item[(M.c)] the speed ${\sigma}^{\diamond}_s$, the left and right 
	states ${u}^{\diamond}_{\rm M}(x^{\diamond}_{s},(n+1/2){\varDelta}t),{u}_{\rm R}(x^{\diamond}_{s},(n+1/2){\varDelta}t)$ satisfy the Rankine--Hugoniot conditions, 
\end{itemize}
where $\bar{u}^{\diamond}_{\rm M}(x)=
\bar{\mathcal U}(x,j{\varDelta}x,u^{\diamond}_{\rm M})$,\;$x^{\diamond}_{p}:=j{\varDelta}x+\sigma^{\diamond}_{p}{\varDelta}
/2$ and $x^{\diamond}_s:=j{\varDelta}x+\sigma^{\diamond}_s{\varDelta}
/2$. 

We denote this approximate Riemann solution, which consists of \eqref{eqn:solution-steady-state}, by ${u}^{\varDelta}(x,t)$. The validity of the above construction is demonstrated in \cite[Appendix A]{T2}.

\begin{remark}\label{rem:middle-time}\normalfont
${u}^{\varDelta}(x,t)$ satisfies the Rankine--Hugoniot conditions
 at the middle time of the cell, $t_{\rm M}:=(n+1/2){\varDelta}t$.
\end{remark}


\begin{remark}\label{rem:approximate}\normalfont
	The approximate solution $u^{\varDelta}(x,t)$ is piecewise smooth in each of the 
	divided parts of the cell. Then, in the divided part, $u^{\varDelta}(x,t)$ satisfies
	\begin{align*}
	(u^{\varDelta})_t+f(u^{\varDelta})_x-g(x,u^{\varDelta})=O(\varDelta x).
	\end{align*}
\end{remark}

\section{Energy inequality}
In this section, we prove Theorem \ref{thm:main2}, i.e., we deduce an energy inequality for our 
solutions in Theorem \ref{thm:main}.
For any fixed $T>0$, we set $N=[T/\varDelta t]$, where $[x]$ is the greatest integer not greater than $x$.  
From \eqref{eqn:cutoff} and finite propagation, we can choose
${R_T}$ large enough such that ${\rm Supp}\:u^{\varDelta}\subset[1,{R_T}]\times[0,T]$. 
Throughout this section, by Landau's symbols such as $O({ \varDelta}x)$,
$O(({ \varDelta}x)^2)$ and $o({ \varDelta}x)$,
we denote quantities whose moduli satisfy a uniform bound depending only on $R_T$ and $M$ in \eqref{eqn:IC}.

From Remark \ref{rem:approximate}, $u^{\varDelta}$ satisfy 
\begin{align*}
\eta_*(u^{\varDelta})_t+q_*(u^{\varDelta})_x-a(x)q_*(u^{\varDelta})=O(\varDelta x)
\end{align*}
on the divided part in the cell where $u^{\varDelta}$ are smooth. Moreover, $u^{\varDelta}$ satisfy an entropy condition (see \cite[Lemma 5.1--Lemma 5.4]{T2}) along 
discontinuous lines approximately. Then, applying the Green formula to $\eta_*(u^{\varDelta})_t+q_*(u^{\varDelta})_x-a(x)q_*(u^{\varDelta})$ in the cell $(j-1){\varDelta}x\leq{x}<(j+1){\varDelta}x,\;n{\varDelta}{t}\leq{t}<(n+1){\varDelta}{t}\quad(j\in J_{n+1},\;n\in{\bf Z}_{\geq0},\;n\leq N)$, we have
\begin{align}
\begin{split}
\eta_*(u^{n+1}_j)\leq&\frac{\eta_*(u^n_{j+1})+\eta_*(u^n_{j-1})}2-\frac{\varDelta t}{2\varDelta x}\left\{q(u^n_{j+1})-q(u^n_{j-1})\right\}\\
&+R(x^n_{j+1},u^n_{j+1}){\varDelta t}-R(x^n_{j-1},u^n_{j-1}){\varDelta t}\\
&+\frac{1}{2\varDelta x}\int^{t_{n+1}}_{t_n}\int^{x_{j+1}}_{x_{j-1}}
\frac{A'(x)}{A(x)}q(u^{\varDelta}(x,t))dxdt+o(\varDelta x),
\end{split}
\label{eqn:energy1}
\end{align}
where $t_n=n{\varDelta t}$ and 
\begin{align*}
R(x,u)=&-\frac{\varDelta x}{4{\varDelta t}}b(x)\left(\frac3{\gamma-1}
(\rho)^{\theta}m+\frac{m^3}{2(\rho)^{\theta+2}} \right)\\
&+\frac{\varDelta t}{4\varDelta x}a(x)\left\{
\frac{\gamma}{\gamma-1}\frac{\rho^{2\theta}m^2}{\rho}+\frac12\frac{m^4}{\rho^3}\right\}\\
&-\frac{\varDelta t}{4\varDelta x}b(x)\left\{\frac{\gamma+\theta+1}{(\gamma-1)\theta}m\rho^{3\theta}+
\frac{\gamma+3\theta+4}{2\theta}\frac{m^3\rho^{\theta}}{\rho^2}+\frac{m^5}{2\rho^{\theta+4}}\right\}.     
\end{align*}

Multiplying the above inequality by $A(x)$, we obtain
\begin{align*}
\sum_{j\in J_{n+1}}\int_{I_j}A(x)\eta_*(u_j
^{n+1})dx\leq A_n+B_n+C_n+o(\varDelta x),
\end{align*}
where $I_j=[(j-1){\varDelta}x,(j+1){\varDelta}x]$ and 
\begin{align*}
A_n=&\sum_{j\in J_{n+1}}\int_{I_j}A(x)\frac{\eta_*(u^n_{j+1})+\eta_*(u^n_{j-1})}2dx
-\sum_{j\in J_{n+1}}\int_{I_j}A(x)\frac{\varDelta t}{2\varDelta x}\left\{q(u^n_{j+1})-q(u^n_{j-1})\right\}dx
,\\
B_n=&-\frac{1}{2\varDelta x}\sum_{j\in J_{n+1}}\int_{I_j}A(x)\left\{\int^{t_{n+1}}_{t_n}\int^{x_{j+1}}_{x_{j-1}}
\frac{A'(y)}{A(y)}q(u^{\varDelta}(y,t))dydt
\right\}
dx,\\
C_n=&\sum_{j\in J_{n+1}}\int_{I_j}A(x)\left\{R(x^n_{j+1},u^n_{j+1})-R(x^n_{j-1},u^n_{j-1})
\right\}{\varDelta t}\:dx.
\end{align*}

We first compute $A_n$.
\begin{align*}
A_n=&\sum_{j\in J_{n+1}}\left\{\int_{I_{j+2}}A(x)dx
+\int_{I_j}A(x)dx\right\}
\frac{\eta_*(u^n_{j+1})}2\\
&+\sum_{j\in J_{n+1}}
\left\{\int_{I_{j+2}}A(x)dx
-\int_{I_j}A(x)dx\right\}
\frac{\varDelta t}{2\varDelta x}q(u_{j+1}^n)\\
=&\sum_{j\in J_n}\left\{\int_{I_{j+1}}A(x)dx
+\int_{I_{j-1}}A(x)dx\right\}
\frac{\eta_*(u^n_{j})}2dx\\
&+\sum_{j\in J_n}
\left\{\int_{I_{j+1}}A(x)dx
-\int_{I_{j-1}}A(x)dx\right\}
\frac{\varDelta t}{2\varDelta x}q(u_{j}^n)\\
=&\sum_{j\in J_n}\int_{I_j}A(x)\eta_*(u^n_{j})dx
+{\varDelta t}\sum_{j\in J_n}\int_{I_j}A'(x)q(u_{j}^n)dx
+O((\varDelta x)^2).
\end{align*}

Next we compute $B_n$. 
Then we have
\begin{equation*}
\begin{split}
B_n=&-\frac{1}{2\varDelta x}\sum_{j\in J_{n+1}}\int_{I_j}A(x)\left\{\int^{t_{n+1}}_{t_n}\int_{I_j}
\frac{A'(y)}{A(y)}q(u^{\varDelta}(y,t))dydt
\right\}dx\\
=&-\sum_{j\in J_{n+1}}
\int^{t_{n+1}}_{t_n}\int_{I_j}A'(x)q(u^{\varDelta}(x,t))dxdt
+O((\varDelta x)^2)\\
=&
-{\varDelta t}\sum_{j\in J_{n+1}}\int_{I_j}A'(x)q(u_{j}^{n+1})dx\\
&+\sum_{j\in J_{n+1}}
\int^{t_{n+1}}_{t_n}\int_{I_j}A'(x)
\left\{q(u^{\varDelta}(x,t_{n+1}+0))-q(u^{\varDelta}(x,t_{n+1}-0))\right\}dxdt\\
&+\sum_{j\in J_{n+1}}
\int^{t_{n+1}}_{t_n}\int_{I_j}A'(x)
\left\{q(u^{\varDelta}(x,t_{n+1}-0))-q(u^{\varDelta}(x,t))   \right\}dxdt+O((\varDelta x)^2)\\
:=&
-{\varDelta t}\sum_{j\in J_{n+1}}\int_{I_j}A'(x)q(u_{j}^{n+1})dx+D_n+E_n+O((\varDelta x)^2).
\end{split}
\end{equation*}Moreover, we find
\begin{align*}
C_n=&-{\varDelta t}\sum_{j\in J_{n+1}}\left\{\int_{I_{j+2}}A(x)dx-\int_{I_j}A(x)dx\right\}R(x^n_{j+1},u^n_{j+1})\\
=&-\varDelta x\:{\varDelta t}\sum_{j\in J_{n+1}}\int_{I_{j+1}}A'(x)dxR(x^n_{j+1},u^n_{j+1})+O((\varDelta x)^2)\\
=&O((\varDelta x)^2).
\end{align*}

Therefore, we have 
\begin{align}
\begin{split}
\sum_{j\in J_{n+1}}\int_{I_j}A(x)\eta_*(u_j
^{n+1})dx
\leq&\sum_{j\in J_{n}}\int_{I_j}A(x)\eta_*(u_j^{n})dx-{\varDelta t}\sum_{j\in J_{n+1}}\int_{I_j}A'(x)q(u_{j}^{n+1})dx\\
&+{\varDelta t}\sum_{j\in J_n}\int_{I_j}A'(x)q(u_{j}^n)dx+D_n+E_n+o(\varDelta x).
\end{split}
\label{eqn:energy2}
\end{align}
Since $\eta_*$ is a convex function, from the Jensen inequality, we obtain 
\begin{align*}
\begin{split}
\sum_{j\in J_{n+1}}\int_{I_j}A(x)\eta_*(u_j
^{n+1})dx
&\leq\sum_{j\in J_{0}}\int_{I_j}A(x)\eta_*(u_j^{0})dx+\sum^n_{k=0}\left(D_n+E_n\right)+o(1)\\
&\leq\int_{\bf R}A(x)\eta_*(u_{0})dx+\sum^n_{k=0}\left(D_n+E_n\right)+o(1).
\end{split}
\end{align*}

Then, we introduce the following proposition:
\begin{proposition}
	\begin{eqnarray}
		\sum^{N-1}_{k=0}\int^{R_T}_1|\overline{u}^{\varDelta}(r,t_k-0)-\overline{u}^{\varDelta}(r,t_k+0)|^2dx\leq C,\label{eqn:Pro1}\\
	\sum^{N}_{k=1}\int^{t_k}_{t_{k-1}}\int^{R_T}_1|\overline{u}^{\varDelta}(r,t_k-0)-\overline{u}^{\varDelta}(r,t)|^2dx=O(\varDelta x)\label{eqn:Pro2}	\end{eqnarray}
\end{proposition}
\eqref{eqn:Pro1} and \eqref{eqn:Pro2} can be obtained in a similar manner to \cite[(6.18)]{T2} and \cite[Lemma 7.1]{T2} respectively.

From the above propositions and the Schwarz inequality, we have
\begin{align*}
\sum^n_{k=0}\left(D_n+E_n\right)=O(\sqrt{\varDelta x}).
\end{align*}Therefore, it follows that
$$
\sum_{j\in J_{n+1}}\int_{I_j}A(x)\eta_*(u_j
^{n+1})dx\leq\int_{\bf R}A(x)\eta_*(u_{0})dx+o(1).
\eqno{(4.5)_{n+1}}
$$

Next, let $s$ be any fixed positive value satisfying $t_{n}\leq s\leq\min\{t_{n+1},T\}$. Applying the Green formula to $\eta_*(u^{\varDelta})_t+q_*(u^{\varDelta})_x-a(x)q_*(u^{\varDelta})$ in the cell $(j-1){\varDelta}x\leq{x}<(j+1){\varDelta}x,\;{t}_n\leq{t}<{s}$, we have
\begin{align*}
\begin{split}
\int_{I_j}\eta_*(u^{\varDelta}(x,s))dx\leq&\frac{\eta_*(u^n_{j+1})+\eta_*(u^n_{j-1})}2-\frac{s-t_n}{2\varDelta x}\left\{q(u^n_{j+1})-q(u^n_{j-1})\right\}\\
&+R(x^n_{j+1},u^n_{j+1})(s-t_n)-R(x^n_{j-1},u^n_{j-1})(s-t_n)\\
&+\frac{1}{2\varDelta x}\int^{s}_{t_n}\int^{x_{j+1}}_{x_{j-1}}
\frac{A'(x)}{A(x)}q(u^{\varDelta}(x,t))dxdt+o(\varDelta x).
\end{split}
\end{align*}
We deduce from the above inequality\setcounter{equation}{5}
\begin{align}
\begin{split}\int_{\bf R}A(x)\eta_*(u^{\varDelta}(x,s))dx=&
\sum_{j\in J_{n+1}}\int_{I_j}A(x)\eta_*(u^{\varDelta}(x,s))dx\\
\leq&\sum_{j\in J_{n}}\int_{I_j}A(x)\eta_*(u_j^{n})dx+O(\varDelta x)
\end{split}\label{eqn:energy3}
\end{align}
in a similar manner to \eqref{eqn:energy2}.

Combining $(4.5)_n$ and \eqref{eqn:energy3}, for $s\leq T$, we conclude
\begin{align}
\begin{split}\int_{\bf R}A(x)\eta_*(u^{\varDelta}(x,s))dx\leq\int_{\bf R}A(x)\eta_*(u_{0})dx+o(1).
\end{split}\label{eqn:energy4}
\end{align}

Then, integrating \eqref{eqn:energy4} over the region $S\in{\bf R}_+$ with $0<\mu(S)<\infty$, we have 
\begin{align}
\begin{split}\int_{\bf R}\int_{S}A(x)\eta_*(u^{\varDelta}(x,s))dxds\leq\mu(S)\int_{\bf R}A(x)\eta_*(u_{0})dx+o(1),
\end{split}\label{eqn:energy5}
\end{align}
where $\mu$ is the one-dimensional Lebesgue measure.

By virtue of the methods of compensated compactness for the approximate solutions 
(see \cite{T8}), there exists a subsequence $u^{\varDelta_k}$ 
such that $(\varDelta x)_k\rightarrow0$ and $u^{\varDelta_k}$ tends to a 
weak solution 
to \eqref{eqn:nozzle} almost everywhere $(x,t)\in{\bf R}_+\times{\bf R}$
as $k\rightarrow\infty$. 
Applying \eqref{eqn:energy5} to the above subsequence and taking the limit, we have 
\begin{align}
\begin{split}\frac1{\mu(S)}\int_{\bf R}\int_{S}A(x)\eta_*(u(x,s))dxds\leq\int_{\bf R}A(x)\eta_*(u_{0})dx.
\end{split}\label{eqn:energy6}
\end{align}
Recalling that $S$ are arbitrary, we have \eqref{eqn:energy}.
Since we can obtain 
\eqref{eqn:energy} for an arbitrary $X$ in \eqref{eqn:cutoff}, we conclude Theorem \ref{thm:main2}.


\appendix

\section{Construction of Approximate Solutions near the vacuum}\label{sec:vacuum}

In this step, we consider the case where $\rho_{\rm M}\leq({\varDelta}x)^{\beta}$,
which means that $u_{\rm M}$ is near the vacuum. In this case, we cannot construct 
approximate solutions in a similar fashion to Subsection 3.1. Therefore, we must
define $u^{\varDelta}(x,t)$ in a different way.

In this appendix, we define our approximate solutions in the cell $(j-1){\varDelta}x\leq{x}<(j+1){\varDelta}x,\;n{\varDelta}{t}\leq{t}<(n+1){\varDelta}{t}\quad(j\in{\bf Z},\;n\in{\bf Z}_{\geq0})$. 
We set $L_j:=-{M}e^{-\int^{(j+1){\varDelta}x}_0b(x)dx}$ and 
$U_j:={M}e^{\int^{j{\varDelta}x}_0b(x)dx}$.

\vspace*{5pt}
{\bf Case 1} A 1-rarefaction wave and a 2-shock arise.

 In this case, we notice that $\rho_{\rm R}\leq ({\varDelta}x)^{\beta},\;
z_{\rm R}\geq L_j$ and $w_{\rm R}\leq U_j$.
\vspace*{5pt}

\begin{center}
{\bf Definition of $\bar{u}^{\varDelta}$ in Case 1}
\end{center}

\vspace*{5pt}
{\bf Case 1.1}
$\rho_{\rm L}>2({\varDelta}x)^{\beta}$

We denote $u^{(1)}_{\rm L}$ a state satisfying $ w(u_{\rm L}^{(1)})=w(u_{\rm L})$ and 
$\rho^{(1)}_{\rm L}=2({\varDelta}x)^{\beta}$. 
Let $u^{(2)}_{\rm L}$ be a state connected to $u_{\rm L}$ on the right by 
$R_1^{\varDelta}(z^{(1)}_{\rm L})(u_{\rm L})$. We set 
\begin{align*}
(z^{(3)}_{\rm L},w^{(3)}_{\rm L})=
\begin{cases}
(z^{(2)}_{\rm L},w(u_{\rm L})),\quad\text{if $z^{(2)}_{\rm L}\geq L_j$},\\
(L_j,w(u_{\rm L})),\quad\text{if $z^{(2)}_{\rm L}< L_j$}.
\end{cases}
\end{align*}

Then, we define 
\begin{align*}
\bar{u}^{\varDelta}(x,t)
=\begin{cases}
R_1^{\varDelta}(z^{(1)}_{\rm L})(u_{\rm L}),\quad
\text{if $(j-1){\varDelta}x
\leq{x}\leq{j}{\varDelta}x+\lambda_1(u^{(2)}_{\rm L})(t-n{\varDelta}t)$}\\\text{and $n{\varDelta}t\leq{t}<(n+1){\varDelta}t$},\vspace*{2ex}\\
\text{a Riemann solution  $(u^{(3)}_{\rm L}$, $u_{\rm R})$},\quad\text{if ${j}{\varDelta}x+\lambda_1(u^{(2)}_{\rm L})(t-n{\varDelta}t)$}\\\text{$<x
\leq (j+1){\varDelta}x$ and $n{\varDelta}t\leq{t}<(n+1){\varDelta}t$}.
\end{cases}
\end{align*}
\begin{figure}[htbp]
\begin{center}
\vspace{0ex}
\hspace{2ex}
\includegraphics[scale=0.38]{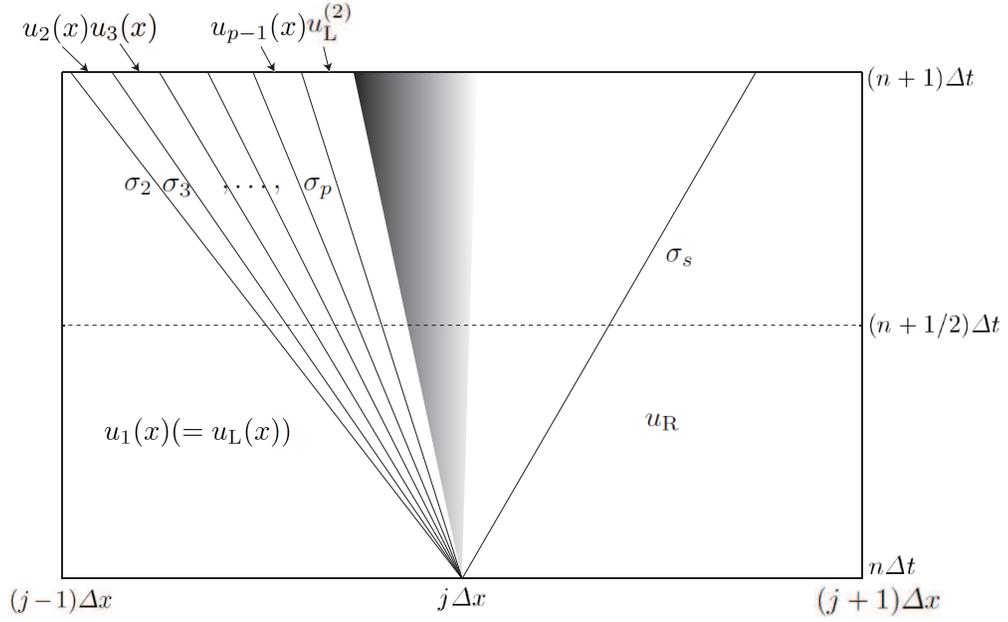}
\end{center}
\caption{{\bf Case 1.1}: The approximate solution $\bar{u}^{\varDelta}$ in the cell.}
\label{Fig:case1.1(ii)}
\end{figure} 

\begin{figure}[htbp]
	\begin{center}
		\vspace{0ex}
		\hspace{2ex}
		\includegraphics[scale=0.32]{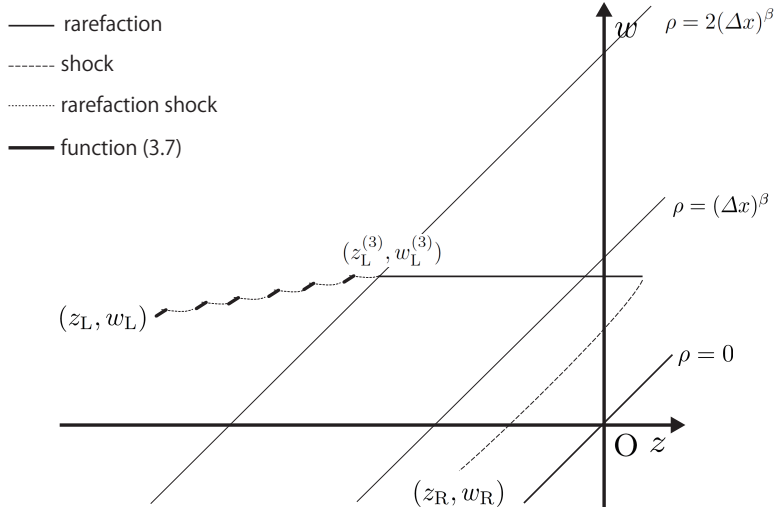}
	\end{center}
	\caption{{\bf Case 1.1}: The approximate solution $\bar{u}^{\varDelta}$ in $(z,w)-$plane.}
	\label{Fig:case1.1(ii)}
\end{figure}

{\bf Case 1.2} $\rho_{\rm L}\leq2({\varDelta}x)^{\beta}$

\vspace*{5pt}
(i) $z(u_{\rm L})\geq{L}_j$\\
In this case, we define $u^{\varDelta}(x,t)$ as a Riemann solution 
$(u_{\rm L},u_{\rm R})$.

\vspace*{5pt}
(ii) $z(u_{\rm L})<L_j$\\
In this case, recalling $z(u_{\rm L})=z(u^n_j)\geq-{M}e^{-\int^{j{\varDelta}x}_0b(x)dx}$, 
we can choose $x^{(4)}$ such that $(j-1){\varDelta}x\leq{x}^{(4)}
\leq(j+1){\varDelta}x$ and 
$
z(u_{\rm L})e^{-\int^{x^{(4)}}_{x_{\rm L}}b(x)dx}=L_j,
$
where $x_{\rm L}:=j{\varDelta}x$.
We set   
\begin{align*}
z^{(4)}_{\rm L}:=z_{\rm L}e^{-\int^{x^{(4)}}_{x_{\rm L}}b(x)dx},\quad{w}^{(4)}_{\rm L}:
=w_{\rm L}e^{-\int^{x^{(4)}}_{x_{\rm L}}b(x)dx}.
\end{align*}

In the region where $(j-1){\varDelta}x
\leq{x}\leq{j}{\varDelta}x+\lambda_1(u^{(4)}_{\rm L})(t-n{\varDelta}t)$ and 
$n{\varDelta}t\leq{t}<(n+1){\varDelta}t$,
we define $\bar{u}^{\varDelta}(x,t)$ as 
\begin{align}
	\bar{z}^{\varDelta}(x,t)=z_{\rm L}e^{-\int^{x}_{x_{\rm L}}b(x)dx},\quad
	\bar{w}^{\varDelta}(x,t)=w_{\rm L}e^{-\int^{x}_{x_{\rm L}}b(x)dx}.
	\label{eqn:vacuum-approximate}
\end{align}

We next solve a Riemann problem
$(u^{(4)}_{\rm L},u_{\rm R})$. In the region where ${j}{\varDelta}x+\lambda_1(u^{(4)}_{\rm L})(t-n{\varDelta}t)\leq{x}\leq(j+1){\varDelta}x$ and 
$n{\varDelta}t\leq{t}<(n+1){\varDelta}t$,  
we define $\bar{u}^{\varDelta}(x,t)$ as this Riemann solution.

We notice that the Riemann solutions in Case 1.2 are also contained in $\varDelta_j$.



\begin{center}
{\bf Definition of ${u}^{\varDelta}$ in Case 1}
\end{center}

\vspace*{5pt}

In the region 
where $\bar{u}^{\varDelta}(x,t)$ is the Riemann solution, we 
define $u^{\varDelta}(x,t)$ by $u^{\varDelta}(x,t)=\bar{u}^{\varDelta}(x,t)$; 
in the region $\bar{u}^{\varDelta}(x,t)$ is \eqref{eqn:vacuum-approximate}, we 
define
\begin{align*}
	\begin{split}
		z^{\varDelta}(x,t)=&\bar{z}^{\varDelta}(x)
		-\left\{a(x)\bar{v}^{\varDelta}(x)(\bar{\rho}^{\varDelta}(x))^{\theta}
		-b(x)\lambda_1(\bar{u}^{\varDelta}(x))\bar{z}^{\varDelta}(x)\right\}
		(t-n{\varDelta{t}}),\\
		w^{\varDelta}(x,t)=&\bar{w}^{\varDelta}(x)+
		\left\{a(x)\bar{v}^{\varDelta}(x)(\bar{\rho}^{\varDelta}(x))^{\theta}
		+b(x)\lambda_2(\bar{u}^{\varDelta}(x))\bar{w}^{\varDelta}(x)\right\}(t-n{\varDelta{t}});
	\end{split}
\end{align*}
otherwise, the definition of 
$u^{\varDelta}(x,t)$ is similar to Subsection 3.1. 
Thus, for a Riemann solution near the vacuum, we define our approximate solution
as the Riemann solution itself.

\vspace*{10pt}
{\bf Case 2} A 1-shock and a 2-rarefaction wave arise.

From symmetry, this case reduces to Case 1.

\vspace*{10pt}
{\bf Case 3} A 1-rarefaction wave and a 2-rarefaction wave arise.

For $u_{\rm L}$ of Case 1, we define $u^*_{\rm L}$ and $\lambda^*_{\rm L}$ as follows. 
\begin{align*}
u^*_{\rm L}=\begin{cases}
u^{(3)}_{\rm L},\quad\text{Case 1.1},\\
u_{\rm L},\quad\text{Case 1.2 (i)},\\
u^{(4)}_{\rm L},\quad\text{Case 1.2 (ii)},
\end{cases}\quad 
\lambda^*_{\rm L}=\begin{cases}
\lambda_1(u^{(2)}_{\rm L}),\quad\text{Case 1.1},\\
\lambda_1(u_{\rm L}),\quad\text{Case 1.2 (i)},\\
\lambda_1(u^{(4)}_{\rm L}),\quad\text{Case 1.2 (ii)}.
\end{cases}
\end{align*}
where $\lambda_1(u)$ be the 1-characteristic speed of $u$. Then, for $u_{\rm L}$ of Case 3, we can determine $u^*_{\rm L}$ and 
$\lambda^*_{\rm L}$ in a similar manner to Case 1. 
From symmetry, for $u_{\rm R}$ of Case 3, we can also 
determine $u^*_{\rm R}$ and $\lambda^*_{\rm R}$.

In the 
region 
$(j-1){\varDelta}x\leq{x}\leq{j}{\varDelta}x+\lambda^*_{\rm L}(t-n{\varDelta}t),\;
{j}{\varDelta}x+\lambda^*_{\rm R}(t-n{\varDelta}t)\leq{x}\leq(j+1){\varDelta}x$ and 
$n{\varDelta}t\leq{t}<(n+1){\varDelta}t$, we define $\bar{u}^{\varDelta}$ in a similar manner to Case 1. In the other 
region, we define $\bar{u}^{\varDelta}$ as the Riemann solution 
$(u^*_{\rm L},u^*_{\rm R})$.

We define ${u}^{\varDelta}$ in the same way as Case 1.

\vspace*{10pt}
{\bf Case 4} A 1-shock and a 2-shock arise.

We notice that $z_{\rm L}\geq L_j,\;w_{\rm L}\leq U_j,\;z_{\rm R}\geq L_j$ and $w_{\rm R}\leq U_j$.
In this case, we define $u^{\varDelta}(x,t)$ as the Riemann
solution $(u_{\rm L},u_{\rm R})$. We notice that the Riemann solution is also contained in $\varDelta_j$.
\vspace*{2ex}

We complete the construction of our approximate solutions.


\end{document}